\documentclass{amsart}

\usepackage{amssymb,amsthm, amsmath}

\title[Characteristic functions and automorphisms of the unit ball]{Characteristic functions for multicontractions and automorphisms of the unit ball}
\author{Chafiq Benhida} 
\address{UFR de Math\'ematiques, Universit\'e des Sciences et Technologies
de Lille, F-59655 Villeneuve D'Ascq Cedex, France}
\email{Chafiq.Benhida@math.univ-lille1.fr}

\author{Dan Timotin}
\address{Institute of Mathematics of the Romanian Academy, P.O. Box 1-764, Bucharest 
014700, Romania}
\email{Dan.Timotin@imar.ro}

\date{\today}

\newtheorem{theorem}{Theorem}[section]
\newtheorem{proposition}[theorem]{Proposition}
\newtheorem{lemma}[theorem]{Lemma}
\newtheorem{corollary}[theorem]{Corollary}

\newtheorem*{thma}{Theorem A}
\newtheorem*{thmb}{Theorem B}

\theoremstyle{definition}
\newtheorem{remark}[theorem]{Remark}
\newtheorem{definition}[theorem]{Definition}

\let\mathbfit\mathbf

\def\H{{\mathbfit H}}

\def\z{{\mathbfit z}}

\def\DD{{\mathcal D}}
\def\EE{{\mathcal E}}
\def\FF{{\mathcal F}}

\def\HH{{\mathcal H}}

\def\KK{{\mathcal K}}
\def\LL{{\mathcal L}}

\def\BBB{{\mathbb B}}
\def\CCC{{\mathbb C}}

\def\HHH{{\mathbb H}}

\def\NNN{{\mathbb N}}

\def\TTT{{\mathbb T}}

\def\0{{\mathbf 0}}
\def\1{{\mathbf 1}}

\def\<{\langle}
\def\>{\rangle}

\def\ker{\mathop{\rm ker}}

\def\Om{\Omega}

\def\Th{{\Theta}}

\def\bflambda{{\pmb\lambda}}

\begin{document}

\begin{abstract}
A \emph{multicontraction} on a Hilbert space $\HH$ is an $n$-tuple
of operators
$T=(T_1,\dots,T_n)$  acting
on $\HH$, such that $\sum_{i=1}^n T_i T_i^*\le \1_\HH$. We obtain some results related to the characteristic function of a commuting multicontraction,
most notably discussing its behaviour with respect to the action of the
analytic automorphisms of the unit ball.
\end{abstract}

\maketitle

\numberwithin{equation}{section}

\section{Introduction}

Let $\HH$ be a Hilbert space; a \emph{multicontraction} is an $n$-tuple
of operators
$T=(T_1,\dots,T_n)$  acting
on $\HH$, such that $\sum_{i=1}^n T_i T_i^*\le \1_\HH$. A theory of dilation
and models for this type of operators has been developed by Gelu Popescu 
in \cite{GP1, GP} and a series of subsequent papers. There is no commutativity 
assumed there, and the isometric dilation obtained is 
related to the Fock space and to representations of the Cuntz
algebra.

Starting mainly with~\cite{Ar},
interest has developed around the case $T$ is formed by commuting operators.
In particular, in the recent paper~\cite{BES}, which is actually
the starting point for this note, a notion of characteristic
function is introduced for commuting tuples, and in a particular
case it is shown that this is a complete unitary invariant.
One computes also, in terms of the characteristic function,
the curvature introduced by Arveson~\cite{Ar2}.
Although some of the results therein follow from the noncommuting
case of Popescu, it is not the case with all of them; moreover,
even when it is, the direct approach might be instructive.

This note investigates further the characteristic function
of a multicontraction. In Section~2 we remind the main definitions
and notations. Section~3 contains some variations around the results in~\cite{BES}.
In Section~4 we investigate a general form of fractional transform,
of which the characteristic function is a particular case.
The main applications are obtained in Section~5, where 
one investigates the relation of the characteristic function
to the automorphisms of the ball applied to a multicontraction.
Here formulas similar to the Moebius transform of a single
contraction are obtained. This is connected to the homogeneous 
operators considered by Misra et al \cite{MS}.

After this paper was completed, we have learnt that further work on closely
related subjects has independently been done by Bhattacharyya, Eschmeier and Sarkar~\cite{ES}.
We will point out, when the case appears, the relation between our results
and~\cite{ES}.

\section{Preliminaries and notations}\label{se:char}

If $\EE_1, \EE_2$ are two Hilbert spaces, and $C:\EE_1\to \EE_2$ is a contraction,
one defines the \emph{defect operator} $D_C=(\1_{\EE_1}-C^*C)^{1/2} \in\LL(\EE_1)$ and the
defect space $\DD_C=\overline{D_C \EE_1}\subset \EE_1$.

Suppose $T=(T_1,\dots, T_n)\in\LL(\HH)^n$ is a commuting multicontraction;
that is, 
\[
\sum_{i=1}^n T_iT_i^*\le 1_\HH. 
\]
This is the same as requiring that 
the row operator $T=(T_1\ \cdots\ T_n):\HH^n\to\HH$ is a contraction. (We will
currently denote with the same letter $T$ the multioperator and the
associated row contraction.) 
Accordingly, we have the operators
$D_T=(\1_{\HH^n}-T^*T)^{1/2}$ and $D_{T^*}=(\1_\HH-TT^*)^{1/2}$, 
and the spaces $\DD_T=\overline{D_T\HH^n}\subset \HH^n$, 
$\DD_{T^*}=\overline{D_{T^*}\HH}\subset\HH$.

For further use, for a multiindex 
$\alpha=(\alpha_1,\dots,\alpha_n)\in \NNN^n$, we shall denote
$|\alpha|=|\alpha_1|+\dots+|\alpha_n|$, and 
$T^\alpha=T_1^{\alpha_1}\cdots T_n^{\alpha_n}$.

If, for $z\in \BBB^n$ (the unit ball of 
$\CCC^n$), the operator $\z:\HH^n\to\HH$ is given by 
\[\z=(z_1 1_\HH\ \cdots\ z_n 1_\HH),\] then $\z$ is a strict contraction,
and thus $1_\HH-\z T^*$ is invertible. We may then define 
\begin{equation}\label{eq:theta0}
\theta_T(z)= -T+D_{T^*}(1_\HH-\z T^*)^{-1}\z D_T:\DD_T\to\DD_{T^*}.
\end{equation}
Thus  $\theta_T(z)$ is an analytic contraction
valued functions defined on $\BBB^n$. 
In \cite{BES}, where $\theta_T(z):\BBB^n\to\LL(\DD_T,\DD_{T^*}) $ is introduced, it is called  the \emph{characteristic function} of $T$, and it is proved (in the commuting case)
that it is a multiplier of the 
corresponding Hardy--Arveson spaces. (For a single contraction all these notions appear in~\cite{SNF}.)

According to a standard terminology introduced in~\cite{SNF} for the
case of a single contraction, we say that two analytic functions
$\Th:\BBB^n\to\LL(\EE_1,\EE_2)$, $\Th':\BBB^n\to\LL(\EE'_1,\EE'_2)$ 
\emph{coincide} if there exist unitary operators $\Omega_i:\EE_i\to\EE'_i$,
$i=1,2$, such that $\Omega_2\Th(z)=\Th'(z)\Omega_1$ for all $z\in\BBB^n$.

The \emph{Hardy-Arveson} space $\H$ is equal to the Hilbert space of analytic
functions on $\BBB^n$ with reproducing kernel $k(z,w)=\frac{1}{1-\< z, w\>}$.
The monomials $z^\alpha$ ($z\in\BBB^n, \alpha\in\NNN^n$) form a complete 
orthogonal family, and we have~\cite[Lemma 3.8]{Ar}
\[
\|z_1^{\alpha_1}\dots z_n^{\alpha_n}\|^2=\frac{\alpha_1!\cdots \alpha_n!}{(\alpha_1+\dots+\alpha_n)!}.
\]
Also, for $\EE$ a Hilbert space, we denote $\H(\EE)=\H\otimes \EE$; thus $\H=\H(\CCC)$. 

The \emph{standard multishift} $S=(S_1,\dots,S_n)$ on $\H$ is defined by $S_i f =z_i f$.
Again in~\cite{Ar} one shows that $S$ is a commuting multicontraction, and $D_{S^*}=P_0$,
where $P_0$ denotes the orthogonal projection onto the constant functions. We will
use the same notation $S$ for the corresponding operators ($S_i\otimes\1_\EE$)
acting on $\H(\EE)$.

If $A\in\LL(\H(\EE),\H(\EE^*))$ is an operator that commutes with the standard
multishift, then $A$ is uniquely defined by its restriction $a$ to $\EE$; we will denote 
then $A=M_a$. One can also view $a$ as a function from $\BBB^n$
to $\LL(\EE,\EE^*)$, by writing $a(z)(\xi)=a(\xi)(z)$.

Now, if $T$ is an arbitrary multicontraction,
we can define a completely positive map 
\begin{equation}\label{eq:rhot}
\rho_T:\LL(\HH)\to \LL(\HH)
\end{equation}
by $\rho_T(X)=\sum_{i=1}^n T_iXT_i^*$, and denote,
as in~\cite{BES}, by $A_\infty(T)\in\LL(\HH)$ the strong limit of the decreasing
sequence of positive operators $\rho_T^k(\1_\HH)$.
The next theorem
appears in~\cite{BES}; most of its ingredients are already present in~\cite{Ar}.

\begin{thma}\label{th:BES}
If $T$ is a commuting contractive tuple of operators on $\HH$, then there exists
a unique bounded linear operator $L:\H(\DD_{T^*})\to\HH$ satisfying 
\begin{equation}\label{eq:defL}
L(f\otimes \xi)=f(T)D_{T^*}\xi
\end{equation}
for all $f$ polynomial, $\xi\in\EE$. 
The adjoint operator $L^*:\HH\to\H(\DD_{T^*})$ is given by
\begin{equation}\label{eq:defL*}
(L^* h)(z)=D_{T^*}(\1_\HH-zT^*)^{-1}h.
\end{equation}
These operators satisfy
\begin{equation}\label{eq:comL}
LS_i=T_iL,\qquad S_i^*L^*=L^* T_i^*.
\end{equation}

We have the identities
\begin{align}
LL^*+A_\infty(T)&=\1_\HH,\label{eq:LA}  \\
L^*L+M_{\theta_T}M_{\theta_T}^*&=\1_{\H(\DD_{T^*})}.\label{eq:Lth}
\end{align}
\end{thma}
When no confusion is possible, we will usually denote simply $A_\infty$
instead of $A_\infty(T)$.

Finally, we remind a dilation result from~\cite{Ar}. A~\emph{spherical operator} 
$Z=(Z_1,\dots, Z_n)$ is a commuting tuple of normal operators such that 
$Z_1Z_1^*+\dots+Z_nZ_n^*=\1$. 

\begin{thmb}\label{th:dilation}
If $T$ is a multicontraction, there is a (essentially uniquely defined) \emph{standard 
minimal dilation}
of~$T$ of the form $S\oplus Z$, where $S$ is a multishift on $\H(\DD_{T^*})$, while
$Z$ is a spherical operator. We have $Z=0$ iff $A_\infty(T)=0$.
\end{thmb}

A few words are in order concerning the relation to the noncommuting
case considered by Popescu in~\cite{GP1, GP} and subsequent papers.
When the contractions $T_i$ do not commute, it is necessary to introduce,
instead of the Hardy--Arveson space $\H(\EE)$, the Fock space
\[
\Gamma(\EE)=\EE \oplus\EE^n \oplus(\EE^n)^{\otimes^2} \oplus\cdots\oplus
(\EE^n)^{\otimes^m}\oplus\cdots. 
\]
One can identify then $\H(\EE)$ with the subspace of $\Gamma(\EE)$
formed by the symmetric tensors; suppose then that $\pi_\EE$ is the orthonormal
projection onto this subspace. The defect spaces of a noncommuting 
multicontraction are defined in a similar manner, and the characteristic function of $T$
introduced by Popescu in~\cite{GP} corresponds to an operator
$\mu_T:\Gamma(\DD_T)\to\Gamma(\DD_{T^*})$ which commutes with the
creation  operators in the Fock space. The relation with the commuting
case is then the formula $M_{\theta_T}=\pi_{\DD_{T^*}} \mu_T|\H(\DD_T)$.
For other connections between the commuting and noncommuting cases, 
one can see~\cite{RBD}, as well as the recent extensive paper~\cite{GP3}.

\section{Classes of multicontractions}


By means of the operator~$A_\infty$ defined above, we can define
some classes of multicontractions, similar to the ones that appear
in~\cite{GP} in the noncommutative case.
 
\begin{definition} 
  The multicontraction~$T$ is called 
\begin{itemize}
\item \em{pure} (or $C_0$) if $A_\infty=0$;
\item $C_1$ if\/ $\ker A_\infty=\{0\}$;
\item completely noncoisometric (c.n.c) if\/ $\ker(\1- A_\infty)=\{0\}$.
\end{itemize}
\end{definition}

The next result is an immediate consequence of Theorem~A.

\begin{proposition}\label{co:cnc}
(i) \cite{BES} $T$ is pure iff $L^*$ is an isometry. In this case 
$T$ is unitarily equivalent to the commuting
tuple $\TTT=(\TTT_1,\dots,\TTT_n)$ on $\HHH_T=\H(\DD_{T^*})\ominus M_{\theta_T}(\H(\DD_T))$
defined by $\TTT_i=P_{\HHH_T}S_i|\HHH_T$.

(ii) $T$ is c.n.c. if and only if $L^*$ is one-to-one. In this case $T^*$ is uniquely
determined by the second equality in~\eqref{eq:comL}.
\end{proposition}

The following lemma has been proved for the noncommuting case in~\cite{GP1} (see 
Remark~2.7 therein).

\begin{lemma}\label{le:kernels}
Suppose~$T$ is a commuting multicontraction.

(i) $\ker A_\infty$ is invariant with respect to $T^*$, and the compression of $T$ to $\ker A_\infty$ is pure.

(ii)
$\ker(\1- A_\infty)$
is invariant with respect to $T^*$ and, if
$\hat T$ is the compression to~$T$ to $\ker(\1_\HH-A_\infty)$, 
then $\DD_{\hat T^*}=\{0\}$.
\end{lemma}

We can  also characterize multishifts by means of their characteristic function.

\begin{proposition}\label{pr:th=0}
If $T$ is pure, then $T$ is unitarily equivalent to the multishift $S$ iff 
$\theta_T\equiv 0$.
\end{proposition}

\begin{proof}

In Theorem~A, if $S$ is the multishift, then $L$ in~\eqref{eq:defL} becomes the 
identity. It follows then from~\eqref{eq:Lth} that $M_{\theta_T}=0$, and thus
$\theta_T\equiv 0$.

Conversely, if $\theta_T\equiv 0$, then from Proposition~\ref{co:cnc} it follows that
$\HHH_T=\H(\DD_{T^*})$, $\TTT_i=S_i$, and $T$ is unitarily equivalent to~$\TTT$.
\end{proof}

\begin{definition}
An analytic function $\Phi:\BBB^n\to\LL(\FF,\EE)$ is called:
\begin{itemize}
\item
\emph{inner} if $M_\Phi:\H(\FF)\to\H(\EE)$
is a partial isometry;
\item
\emph{outer} if $M_\Phi:\H(\FF)\to\H(\EE)$
has dense range.
\end{itemize}
\end{definition}

The following characterizations of classes of commuting multicontractions,
by means of their characteristic functions, are
similar to those obtained in~\cite{SNF} for one contraction and in~\cite{GP} for
noncommuting contractions.

\begin{theorem}\label{th:classes}
Suppose~$T$ is a c.n.c. multicontraction on~$\HH$. Then:
\begin{enumerate}
\item $T$ is pure iff $\theta_T$ is inner.
\item $T$ is of class $C_1$ iff $\theta_T$ is outer.
\end{enumerate}
\end{theorem}

\begin{proof}
(1) If $T$ is pure, then $L$ is a coisometry, and thus $L^*L$ is a projection.
From~\eqref{eq:Lth} it follows that $M_{\theta_T}M_{\theta_T}^*$ is also a
projection; thus $M_{\theta_T}$ is a partial isometry and $\theta_T$ is inner.

Conversely, if $M_{\theta_T}$ is a partial isometry, then~\eqref{eq:Lth} implies
that $L$ is a partial isometry, and then from~\eqref{eq:LA} it follows
that $A_\infty$ is a projection. Since $T$ is c.n.c., we must have 
$\ker(\1- A_\infty)=\{0\}$.
Thus $A_\infty=0$, which means that $T$ is pure.

(2) Note that $\ker A_\infty=\{0\}$, is equivalent, by~\eqref{eq:LA}, to
$\ker(\1-LL^*)=\{0\}$. But it is easy to see (for any bounded operator
$L$, actually) that this last relation is equivalent to $\ker(\1-L^*L)=\{0\}$.
By~\eqref{eq:Lth}, this is the same as $\ker(M^*_{\theta_T})=\{0\}$, or 
$\theta_T$ outer.
\end{proof}

To end the section, let us note that the model provided by Proposition~\ref{co:cnc}~(i)
for pure contractions 
can be extended up to a certain point to c.n.c. multicontractions, in a manner similar to~\cite{SNF}
(for a single contraction) or to~\cite{GP} (for noncommuting multicontractions).  We give only some indications in this direction, 
mostly in order to obtain an extension of Proposition~\ref{pr:th=0}. More details, including
an investigation of the model
space, can be found in~\cite{ES}.

Remember that $\DD_{M_{\theta_T}}=\overline{(I-M_{\theta_T}^*M_{\theta_T})^{1/2}\H(\DD_T)} \subset\H(\DD_T)$.
Consider the space 
$\KK_T= \H(\DD_{T^*})\oplus \DD_{M_{\theta_T}}$
and the two mappings $v:\H(\DD_T)\to\KK_T$, $v_*:\H(\DD_{T^*})\to\KK_T$, 
defined by
\begin{equation}\label{eq:vv*}
v(f)=M_{\theta_T}f\oplus D_{M_{\theta_T}}f,\qquad v_*(f)=f\oplus 0.
\end{equation}
It is easy to check that $v,v_*$ are isometries, that
\begin{equation}\label{eq:span2}
\KK_T=v(\H(\DD_T)) \vee v_*(\H(\DD_{T^*})),
\end{equation}
and that $v_*^* v=M_{\theta_T}$; consequently
$v^*v_*=M_{\theta_T}^*$.
Define 
$\HHH_T= \KK_T\ominus u(\H(\DD_T))$.

If $k\in\HHH_T$, and $k\perp v_*(\H(\DD_{T^*}))$, we must have $k=0$ 
by~\eqref{eq:span2}. Thus the projection $P_{\HHH_T}$ onto $\HHH_T$
has dense range when restricted to $v_*(\H(\DD_{T^*}))$.
Also, if $f\in\H(\DD_{T^*})$, then
\[
\|v_* f\|^2= \|P_{\HHH_T} v_* f\|^2+\|v^*v_* f\|^2=
\|P_{\HHH_T} v_* f\|^2+ \|M_{\theta_T}^* f\|^2.
\]

On the other hand, by~\eqref{eq:Lth}
for any $f\in\H(\DD_{T^*})$ we have
$\|Lf\|^2+\|M_{\theta_T}^* f\|^2=\|f\|^2$.
Therefore, the map $ Lf\mapsto P_{\HHH_T} v_* f$ is an isometry. We
have just noticed that its range is dense in $\HH$; but its domain is also
dense in $\HHH_T$ by Corollary~\ref{co:cnc}. We obtain then a unitary 
$\Phi:\HH\to\HHH_T$, defined by the formula 
\begin{equation}\label{eq:defPhi}
\Phi(Lf)=P_{\HHH_T} v_* f.
\end{equation}

Now, since $\1_\KK- P_{\HHH_T}=vv^*$, we have, using~\eqref{eq:Lth},
\begin{equation*}
\begin{split}
v_*^* (\Phi Lf)&=
v_*^* P_{\HHH_T} v_*=\1_{\H(\DD_{T^*})} - v_*^*(\1_\KK- P_{\HHH_T}) v_*\\
&=\1_{\H(\DD_{T^*})}- v_*^*vv^* v_*
=\1_{\H(\DD_{T^*})}- M_{\theta_T}M_{\theta_T}^*= L^*L.
\end{split}
\end{equation*}
Again, since the range of $L$ is dense in $\HH$, it follows that
\begin{equation}\label{eq:11}
	v_*^*\Phi=L^*.
\end{equation}

We may then define a multioperator $\TTT$ on $\HHH_T$ 
by requiring that 
$v_*^* \TTT^*_i k=S_i^* v_*^* k$.
Applying~\eqref{eq:11} and~\eqref{eq:comL}, we have
\[
v_*^*\TTT_i^*\Phi h=S_i^* v_*^* \Phi h=S_i^* L^* h= L^* T_i^* h
=v_*^* \Phi T_i^* h;
\]
since $v_*^*$ is one-to-one, this 
shows that $\TTT$ is a multicontraction unitarily equivalent to $T$.
Using this unitary equivalence, one can prove, on
the lines
of Theorem 4.4 in~\cite{BES}, that the characteristic
function is a complete unitary invariant for c.n.c. contractions.

\begin{theorem}\label{th:genmodel}
Two c.n.c. contractions
are unitarily equivalent if and only if their characteristic functions
coincide.
\end{theorem}

By Lemma~\ref{le:kernels} (ii),
this is a natural framework
for the extension of~\cite[Theorem 4.4]{BES}. Note also that an alternate
proof of Theorem~\ref{th:genmodel} can be obtained by using the noncommutative
theory of~\cite{GP}.

As a consequence, 
we can obtain a generalization of Proposition~\ref{pr:th=0}.

\begin{corollary}\label{co:th0}
If $T$ is c.n.c., then $T$ is unitarily equivalent to the multishift $S$ iff 
$\theta_T\equiv 0$.
\end{corollary}

\begin{remark}
The main drawback of Theorem~\ref{th:genmodel}
is that not all contractive multipliers coincide with characteristic 
functions of commuting multicontractions, and, contrary to the noncommuting
case of~\cite{GP}, we do not know of a simple way to characterize
those that do.
 A simple example is given by the 
null characteristic function: $\theta:\BBB^n\to\LL(\EE_1,\EE_2)$ defined
by $\theta(z)=0$ for all $z$ coincides with a characteristic function 
if and only if $\dim\EE_1=\infty$. Indeed, it is obvious that in this
case coincidence is equivalent to equalities of the dimensions of the 
domain  and of the range. On the other hand, Corollary~\ref{co:th0}
implies that, if $\theta$ coincides with $\theta_T$, then $T$ has to be
the multishift $S$ on some space $\H(\EE)$. But then $\dim\DD_{S^*}=\dim\EE$
(and is thus arbitrary), while (for $n\ge 2$) $\dim\DD_S=\infty$. This follows immediately
from the fact that $\dim\ker\DD_S=\infty$, since it contains, for instance,
all elements in $\bigoplus_{i=1}^n \H(\EE)$ 
of the form $ (z_2 f)\oplus (-z_1 f)\oplus\ \bigoplus_{i=3}^n 0$ (with $f\in\H(\EE)$).
\end{remark}

\section{Fractional transforms}

It is  useful to regard characteristic functions of multicontractions in a 
larger context, namely as a particular case of fractional transforms.
This section is a development of some results in~\cite{BT2}.

Let $A,W\in\LL(\EE_1,\EE_2)$  
be contractions such that the inverse
$(I+WA^*)^{-1}$ exists;
define the operator
$\Psi_A(W)\in\LL(\EE_1,\EE_2 )$ by the formula
\begin{equation}\label{eq:5}
\Psi_A(W)=A+D_{A^*}(I+WA^*)^{-1}WD_A
\end{equation}
The inverse exists if for example  $\|WA^*\|<1$.
A related operator is
\begin{equation}\label{eq:5bis}
\psi_A(W)=\Psi_A(W)|{\DD_{A}}:\DD_A \to \DD_{A^*};
\end{equation}
note also that $\Psi_A(W)|\DD_A^\perp=A|\DD_A^\perp$, and it maps this
subspace unitarily onto $\DD_{A^*}$.

We will use repeatedly the formulas
\begin{equation}\label{eq:minus}
\Psi_{-A}(-W)=-\Psi_A(W), \quad \psi_{-A}(-W)=-\psi_A(W).
\end{equation}

\begin{lemma}\label{le:2.1}
 {\rm (i)} We have the relations
\begin{align*}
I-\Psi_A(W)^*\Psi_A(W)&=
D_A(I+W^*A)^{-1}(I-W^*W)(I+A^*W)^{-1}D_A,\\
I-\Psi_A(W)\Psi_A(W)^*&=
D_{A^*}(I+WA^*)^{-1}(I-WW^*)(I+AW^*)^{-1}D_{A^*}.
\end{align*}
In particular, $\|\Psi_A(W)\|\le 1$, and, if $W$ is an isometry
(or coisometry), then $\Psi_A(W)$ and $\psi_A(W)$ are isometries (or coisometries,
respectively).
\end{lemma}

The lemma is proved by straight computation. As a consequence, we can define 
isometric operators 
$\Omega:\DD_{\Psi_A(W)}\to\DD_W$ and
$\Omega_*:\DD_{\Psi_A(W){}^*}\to\DD_{W^*}$ by
\begin{align}
\Omega{D_{\Psi_A(W)}}x & =D_W(I+A^*W)^{-1}{D_A}x,\label{eq:omega} \\
\Omega_*{D_{\Psi_A(W){}^*}}x & ={D_{W^*}}(I+AW^*)^{-1}{D_{A^*}}x.\label{eq:omega*}
\end{align}

\begin{remark}\label{re:w0}
There is an important case when we can strengthen these statements.  Suppose
that $W_0:\DD_A\to\DD_{A^*}$ and $I+W_0A^*:\DD_{A^*}\to \DD_{A^*}$ is invertible. (Note that
$A(\DD_A)\subset \DD_{A^*}$ and $A^*( \DD_{A^*})\subset\DD_A$.) Then, if 
$W\in\LL(\EE_1,\EE_2)$ is defined by $W=W_0 P_{\DD_A}+A P_{\DD_A^\perp}$, then $I+WA^*$ is invertible,
$\DD_W\subset \DD_A$, $\DD_{W^*}\subset\DD_{A^*}$, and $\Omega,\Omega_*$ are actually unitary.
In this case formulas~\eqref{eq:omega} and~\eqref{eq:omega*} yield identifications of the 
defect spaces of $\Psi_A(W)$.
\end{remark}

\begin{proposition}\label{pr:2.1}
Suppose the operators $I+VA^*$, $I+WA^*$, $I+\Psi_A(V)W^*$, $I+V\Psi_A(W)^*$ are all 
invertible. Then, if $\Om,\Om_*$ are defined by~\eqref{eq:omega}
and~\eqref{eq:omega*}, we have
\begin{equation}\label{eq:2.4}
\Om_*\psi_{\Psi_A(W)}(V)= \psi_W({\Psi_A(V)})\Om.
\end{equation}
\end{proposition}

It should be noted that $\Om$ and $\Om_*$ depend only on $A$ and $W$ (and not on~$V$).

\begin{proof}
%
%
Since the two terms of~\eqref{eq:2.4} act on 
$\DD_{\Psi_A(W)}$,  in order to check it we have to apply
them to  ${D_{\Psi_A(W)}}x$. 
The proof is a rather tedious computation for which we give only
some indications.

One checks first the two formulas
\begin{equation*}
\Psi_A(W) D_{A}=D_{A{}^*}(I+W A{}^*)^{-1}(W+A)
\end{equation*}
and
\begin{equation*}
\begin{split}
&D_{W^*}(I+\Psi_A(V)W^*)^{-1}D_{A^*}(I+V A^*)^{-1}\\
&\qquad= D_{W^*} (I
+AW^*)^{-1}D_{A^*}(I+V{\Psi_A(W)}^*)^{-1}.
\end{split}
\end{equation*}
Using them, one shows that
\begin{equation}\label{eq:2.1}
\begin{split}
&\Omega_*{\psi_{\Psi_A(W)}}(V){D_{\Psi_A(W)}}x\\
&\quad={D_{W^*}}(I+\Psi_A(V)W^*)^{-1}\big(\Psi_A(V) D_{A}-D_{A^*}(I+V
A^*)^{-1}
 D_{A^*}W(I+A^*W)^{-1}D_A\big).
\end{split}
\end{equation}

On the other hand, writing explicitely the left term in~\eqref{eq:2.4} yields
\begin{equation*}
\begin{split}
&{\psi_W ({\Psi_A(V)})}{\Omega D_{\Psi_A(W)}}x\\
&\quad=WD_W(I+A^*W)^{-1}D_Ax
+D_{W^*}(I+\Psi_A(V)W^*)^{-1}\Psi_A(V)(I+A^*W)^{-1}D_Ax \\
&\qquad\qquad-D_{W^*}(I+\Psi_A(V)W^*)^{-1}\Psi_A(V)W^*W(I+A^*W)^{-1}D_
Ax. 
 \end{split}
\end{equation*}
Since
\begin{equation*}
\begin{split}
&D_{W^*}(I+\Psi_A(V)W^*)^{-1}\Psi_A(V)W^*W(I+A^*W)^{-1}D_Ax\\
&\quad=-D_{W^*}(I+\Psi_A(V)W^*)^{-1}W(I+A^*W)^{-1}D_Ax +
WD_W(I+A^*W)^{-1}D_Ax,
\end{split}
\end{equation*}
it follows that
\begin{equation*}
\begin{split}
&{\psi_W ({\Psi_A(V)})}{\Omega D_{\Psi_A(W)}}x
\\
&\quad=D_{W^*}(I+\Psi_A(V)W^*)^{-1}[\Psi_A(V)(I+A^*W)^{-1}D_Ax+W(I+A^*W)^{-1}D_Ax
]\\
&\quad=D_{W^*}(I+\Psi_A(V)W^*)^{-1}[\Psi_A(V) D_{A}x+(-\Psi_A(V)A^* + I)W(I+A^*W)^{-1}D_Ax]\\
&\quad=D_{W^*}(I+\Psi_A(V)W^*)^{-1}[\Psi_A(V)
D_{A}-D_{A^*}(I+V A^*)^{-1}
 D_{A^*}W(I+A^*W)^{-1}D_A ]x.
\end{split}
\end{equation*}
Comparing this last equality with~\eqref{eq:2.1}, one obtains the desired equality.
\end{proof}

The relation with the characteristic function is obtained by noting
first that, for $z\in \BBB^n$, the operator $\z:\HH^n\to\HH$ 
is a strict contraction,
and thus $1_\HH-\z T^*$ is invertible. Then, by comparing~\eqref{eq:5bis}
and~\eqref{eq:theta0}, we see that
\begin{equation}\label{eq:theta}
\theta_T(z)=\psi_{-T}(\z).
\end{equation}

\section{Involutive automorphisms of the unit ball}

A main interest of formula~\eqref{eq:theta} is that it allows 
to make the connection between characteristic functions and 
automorphisms of the ball applied to multicontractions.

The involutive automorphisms of the unit ball $\BBB^n$ are defined~\cite{Ru} by
\[
\phi_\lambda(z)=\lambda-\frac{s_\lambda}{1-\<z,\lambda\>}
(z-(1-s_\lambda)P_\lambda z),
\]
where $\lambda \in\BBB^n$, $s_\lambda =(1-|\lambda |^2)^{1/2}$, and $P_\lambda $ is the projection onto
the space spanned by~$\lambda $.
Among the properties of these maps, we note that $\phi_\lambda$ is involutive,
that is, $\phi_\lambda\circ\phi_\lambda = 1_{\BBB^n}$, that it maps the
unit ball onto the unit ball, and the unit sphere onto the unit sphere.

The relation with the previous section  is given by the next proposition,
whose proof is a direct computation.

\begin{proposition}\label{pr:alpha}
If we identify  $\lambda ,z\in\BBB^n$ with strict contractions in
$\LL(\CCC^n,\CCC)$, then 
\[
\phi_\lambda (z)=\Psi_{\lambda} (-z)=\psi_{\lambda} (-z).
\]
\end{proposition}

Denote $\bflambda=\lambda\otimes 1_\HH:\HH^n\to\HH$; then $\bflambda$
is a strict contraction. 
If $T=(T_1,\dots,T_n)$ is a multicontraction, then $\bflambda T^*$ is a strict
contraction, and we may define $\Psi_{\bflambda}(-T)=\psi_{\bflambda}(-T)$.

In case $T$ is commutative, $\phi_\lambda(T)$ is
defined by the analytic functional calculus, and it is in turn a commuting
row contraction. Moreover
\[
\phi_\lambda (T)=\Psi_\bflambda(-T).
\]

There are several properties of the multicontraction $T$ that
are also inherited by $\phi_\lambda (T)$. The first one can
be proved directly.

\begin{proposition}\label{pr:phicnc}
With the above notations, $\ker(1-A_\infty(T)) =\ker(1-A_\infty(\phi_\lambda (T)))$.
In particular, $T$ is c.n.c., if and only if $\phi_\lambda (T)$ is c.n.c.
\end{proposition}

\begin{proof} Denote $\phi_\lambda (T)=R=(R_1,\dots,R_n)$.
Since all $R_i$ are analytic functions in
$T_1,\dots,T_n$, all subspaces of $\HH$ invariant to $T$ are also
invariant to $R$. But then the relation
$\phi_\lambda(R)=T$ implies that $R$ and $T$ have the same invariant
subspaces, and the same is true for $T^*$ and $R^*$.

Suppose then that $\KK=\ker(1-A_\infty(T))$ is not contained in $\ker(1-A_\infty(R))$.
Since the elements of this last subspace are characterized by the
fact that $\rho_R^k(\1)=\1$ for all $k$ ($\rho_R$ as defined by~\eqref{eq:rhot}),
there exists $x\in\ker(1-A_\infty(T))$ and a first index $K\ge 1$ for which
$\rho_R^K(\1)(x)\not=x$. Since $\rho_R^k(\1)\le \1$ for all $k$, and $\KK$
is invariant to $R^*$, it follows that we can find $y\in\KK$, $y\not=0$ (of the
form $y=R_{i_1}^*\dots R_{i_{K-1}}^*x$), such that $\sum_{i=1}^n\|R_i^* y\|^2<\|y\|^2$.
On the other side, $y\in\KK$ implies $\sum_{i=1}^n\|T_i^* y\|^2=\|y\|^2$.

Consider now the multioperators $T',R'$, which are the compressions 
of $T$ and $R$ respectively to $\KK$.
It is easy to see
that $R'=\phi_\lambda (T')=\Psi_\bflambda(-T')$. But the definition
of $\KK$ implies that $T'$, as a row contraction operator, is a coisometry.
By Lemma~\ref{le:2.1}, $R'$ should also be a coisometry, which contradicts
the existence of $y$. 

We have thus proved that $\ker(1-A_\infty(T)) \subset\ker(1-A_\infty(R))$.
But then $\phi_\lambda(R)=T$ implies that we actually have equality.
\end{proof}

According to~\eqref{eq:omega} and \eqref{eq:omega*}, we have operators $\Om:\DD_{\phi_\lambda (T)}\to \DD_T$ and 
$\Om_*:\DD_{\phi_\lambda (T)^*}\to \DD_{T^*}$; moreover,
since Remark~\ref{re:w0} applies,   $\Om,\Om_*$ 
are actually unitary maps. They provide identifications of the 
defect spaces of $ \phi_\lambda (T)$; what is more interesting,
they can be used in order to obtain a formula for the characteristic function of $\phi_\lambda(T)$.

\begin{theorem}\label{th:theta} (i) The operators $\Om$ and $\Om_*$ are unitaries. 

(ii) $\theta_{\phi_\lambda(T)}$ coincides with
$\theta_T(\phi_{\lambda} (z))$.
\end{theorem}

\begin{proof}
(i) follows from the fact that we can apply Remark~\ref{re:w0} ($A=\mathbf\lambda$ is
a strict contraction). As for (ii),
by Proposition~\ref{pr:alpha}, formulas~\eqref{eq:theta} and~\eqref{eq:minus}, we have
\[
\theta_{\phi_\lambda(T)}(z)= \theta_{\Psi_\bflambda(-T)}(z)
=\psi_{-\Psi_\bflambda(-T)}(\z)= -\psi_{\Psi_\bflambda(-T)}(-\z).
\]
But Proposition~\ref{pr:2.1} applied to the case $A=\bflambda$, 
$V=-\z$, $W=-T$ says that $\psi_{\Psi_\bflambda(-T)}(-\z)$ 
(and thus also $-\psi_{\Psi_\bflambda(-T)}(-\z)$) coincides
with $\psi_{-T}(\Psi_\bflambda(-\z))$. Note that the unitaries
$\Omega, \Omega_*$ in Proposition~\ref{pr:2.1} do not depend on 
$V$, and thus in our case the unitaries implementing the 
coincidence do not depend on~$z$.

Finally, applying again
Proposition~\ref{pr:alpha} and formula~\eqref{eq:theta}, we obtain
\[
\psi_{-T}(\Psi_\bflambda(-\z))= \psi_{-T}(\phi_\lambda(z) \otimes \1_\HH)
= \theta_T (\phi_\lambda(z)),
\]
which ends the proof.
\end{proof}

Theorem~\ref{th:theta} is the generalization of the well known
formula for the characteristic function of a Moebius transform
of a single contraction~\cite[VI.1.3]{SNF}. However,
the definition of the automorphism of the ball makes
$\phi_\lambda$ involutive, and thus $\phi_\lambda^{-1}=\phi_\lambda$.
The change of sign in the usual definition of the Moebius transforms
accounts for the apparition in~\cite{SNF} of a slightly different
formula. 

Note also that a weaker result along the lines of Theorem~\ref{th:theta}
appears in~\cite{ES}.

As a first application, we obtain a partial extension of the relation
between the spectrum and the characteristic function that exists 
for single
contractions~\cite{SNF}.
Recall (see, for instance,~\cite{Mu}) that, for a commuting multioperator $T=(T_1,\dots, T_n)$, one can define
its \emph{right spectrum} by:
\[
\sigma_r(T)=\{\lambda\in\CCC^n : \sum_{i=1}^n (T_i-\lambda_i)(T_i-\lambda_i)^*
\text{ is not invertible}\}.
\]

\begin{proposition}
If $\lambda\in\BBB^n$, then $\lambda\in\sigma_r(T)$ iff $\theta_T(-\lambda)$ is not surjective.
\end{proposition}

\begin{proof} Note that, since $\sum_{i=1}^n T_i T_i^*$ is invertible if
and only if $(T_1 \cdots T_n)$ is surjective, we have 
\[
\sigma_r(T)=\{\lambda\in\CCC^n : \big( (T_1-\lambda_1) \cdots
(T_n-\lambda_n)\big)
\text{ not surjective}\}.
\]
Since $T$ maps unitarily $\DD_T^\perp$ onto $\DD_{T^*}^\perp$, it
follows that $0\in\sigma_r(T)$ iff $\theta_T(0)$ is not surjective.
Thus the claim is true for $\lambda=0$. For other values of $\lambda$, since, by the
spectral mapping theorem,  
\[
\sigma_r(\phi_\lambda(T))=\phi_\lambda(\sigma_r(T)),
\]
we have $\lambda\in\sigma_r(T)$ iff $0\in\sigma_r(\phi_\lambda(T))$
(note that $\phi_\lambda(\lambda)=0$). This is
equivalent to $\theta_{\phi_\lambda(T)}(0)$ not surjective. Since, by 
Theorem~\ref{th:theta}, $\theta_{\phi_\lambda(T)}(0)$ is unitarily
equivalent to $\theta_T(\phi_{-\lambda}(0))$, and $\phi_{-\lambda}(0)=-\lambda$,
the proposition is proved.
\end{proof}

\begin{remark}
Naturally, a corresponding result can be proved for the \emph{left spectrum}
$\sigma_l(T):=\overline{\sigma_r(T^\sharp)}$, where $T^\sharp=(T_1^*,\dots,T_n^*)$. This would however require the 
assumption that $T^\sharp$ is a multicontraction. If one would want to deduce
a consequence about the \emph{Harte spectrum} $\sigma_H(T)=\sigma_r(T)
\cup \sigma_l(T)$, one should assume \emph{both} $T$ and $T^\sharp$ 
multicontractions, which is a rather unnatural hypothesis (for instance,
it is not satisfied by the multishift for $n\ge 2$). 
\end{remark}

The next consequence concerns the multishift.

\begin{proposition}\label{pr:alphaS}
If $S$ is a multishift and $\lambda\in\BBB^n$, then $\phi_\lambda(S)$ is also
a multishift, of the same multiplicity.
\end{proposition}

\begin{proof}
By Proposition~\ref{pr:phicnc} $\phi_\lambda(S)$ is a c.n.c. multicontraction, 
while Theorem~\ref{th:theta} implies that its
characteristic function is identically zero. 
We may then apply Corollary~\ref{co:th0} to conclude that $\phi_\lambda(S)$
is a multishift. The equality of the multiplicities follows from the equality
of the defect spaces of $S$ and $\phi_\lambda(S)$, as given by Theorem~\ref{th:theta},~(i).
\end{proof}

We can also study the relation with
the model spaces. 

\begin{lemma}\label{le:1}
(i) If $Z$ is spherical, then $\phi_\lambda(Z)$
is also spherical for all $\lambda\in\BBB^n$.

(ii) If the multicontraction $V$ is a minimal dilation for $T$, then 
$\phi_\lambda(V)$ is a minimal dilation for $\phi_\lambda(T)$.
\end{lemma}

\begin{proof}
(i) follows immediately from the functional calculus  for commuting normal
operators. As for (ii), one sees easily that 
$\phi_\lambda(V)$ is a dilation for $\phi_\lambda(T)$, and that
the space spanned by the
powers of $\phi_\lambda(V)$ applied to $\HH$ is contained in the one
spanned by powers of $V$. On the other side, $\phi_\lambda$
is involutive, which gives the opposite relation.
\end{proof}

\begin{proposition}\label{pr:classes}
(i) If the standard minimal dilation of $T$ is $S\oplus Z$, then the standard 
minimal dilation of $\phi_\lambda(T)$ is $\phi_\lambda(S)\oplus\phi_\lambda( Z)$.

(ii) If $T$ is pure, then $\phi_\lambda(T)$ is pure.

(iii) If $T$ is of class $C_1$, then $\phi_\lambda(T)$ is of class $C_1$.
\end{proposition}

\begin{proof}
(i) and (ii) follow immediately from Proposition~\ref{pr:alphaS} and
Lemma~\ref{le:1}. As for (iii),
we will apply Theorem~\ref{th:classes}. If $T$ is of class $C_1$, then $\theta_T$
is outer. To show that $\theta_{\phi_\lambda(T)}$ is also outer, note
first that by Theorem~\ref{th:theta} it is enough to show that
$\theta_T(\phi_{\lambda} (z))$ is outer. Denote by $C_\lambda$ the operator
$f\mapsto f(\phi_{\lambda} (z))$; it is invertible, since $C_\lambda^2=\1$. Then
\[
M_{\theta_T(\phi_{\lambda} (z))}C_\lambda=C_\lambda M_{\theta_T},
\]
whence it follows that, if $M_{\theta_T}$ has dense range, then 
$M_{\theta_T(\phi_{\lambda} (z))}$ also has.
\end{proof}

\section{General automorphisms}

The general form of an automorphism $\alpha$ of $\BBB^n$ is 
\[
\alpha=\omega\circ \phi_{\lambda},
\]
for $\lambda\in\BBB^n$, and $\omega$ a unitary map of $\CCC^n$
(see, for instance,~\cite[Theorem~2.2.5]{Ru}). We may then complete
the results of the previous section by taking into account the
action of the unitary $\omega$; since we  identify
$\CCC^n$ with row $1\times n$ matrices, and regard $T$ as a row operator,
it is natural to consider the action of $\omega$ as matrix multiplication to
the right.

\begin{lemma}\label{le:omega} Suppose $T'=T ( \1_\HH\otimes\omega)$. Then
$\rho_{T'}=\rho_T$, while $\theta_{T'}(z)$ coincides with $\theta_T( z\omega^*)$.
\end{lemma}

\begin{proof}
We can write $\rho_T(X)= T (X\otimes \1_{\CCC^n} )T^*$; then
\[
\rho_{T'}(X)= T(\1_\HH\otimes\omega) (X\otimes \1_{\CCC^n})(1_\HH\otimes\omega^*) T^*=
T (X\otimes \1_{\CCC^n}) T^*=\rho_T(X).
\]

As concerns the defect spaces and operators, 
we have $D_{T'{}^*}=D_{T^*}$ and $\DD_{T'{}^*}=\DD_{T^*}$, while
$D_{T'}=(1_\HH\otimes\omega^*) D_T (\1_\HH\otimes\omega)$ and
$\DD_{T'}=(1_\HH\otimes\omega^*) \DD_T$. Consequently
\begin{align*}
\theta_{T'}(z)&= -T'+D_{T'{}^*}(1_\HH-\z T'{}^*)^{-1}\z D_{T'}\\
&=-T(1_\HH\otimes\omega)+D_{T^*}(1_\HH-\z(1_\HH\otimes\omega^*) T^*)^{-1}
\z(1_\HH\otimes\omega^*) D_T(1_\HH\otimes\omega)\\
&=\theta_T(z\omega^*) (1_\HH\otimes\omega),
\end{align*}
and the lemma is proved. 
\end{proof}

\begin{corollary}\label{co:omega}
Suppose $T'=T ( \1_\HH\otimes\omega)$. Then:

(i) $T'$ is c.n.c. ($C_1, C_0$) iff $T$ is c.n.c. ($C_0, C_1$, respectively).

(ii) If $T$ is a multishift, then $T'$ is a multishift of the same multiplicity.
\end{corollary}

\begin{proof}
The results in (i) are immediate consequences of the equality $\rho_{T'}=\rho_T$,
while for (ii) we have to use, besides Lemma~\ref{le:omega},
Corollary~\ref{co:th0}.
\end{proof}

Gathering the results in Propositions~\ref{pr:phicnc},~\ref{pr:alphaS},~\ref{pr:classes},
Lemma~\ref{le:omega}, Corollary~\ref{co:omega},
and Theorem~\ref{th:theta}, we obtain a general result concerning the
action of an analytic automorphism of the unit ball on a multicontraction.

\begin{theorem}\label{th:last}
Suppose $\alpha:\BBB^n\to\BBB^n$ is an analytic automorphism, while $T$
is a multicontraction. Then:

(i) $\alpha(T)$ is c.n.c. ($C_1, C_0$) iff $T$ is c.n.c. ($C_0, C_1$, respectively).

(ii) $\theta_{\alpha(T)}$ coincides with $\theta_T\circ\alpha^{-1}$.

(iii) If $T$ is a multishift, then $\alpha(T)$ is a multishift of the same multiplicity.
\end{theorem}


In~\cite{MS} the notion of \emph{homogeneous} $n$-tuples of operators
is introduced in a general context. In our case, a multicontraction
$T$ is homogeneous if $\alpha(T)$ is unitarily equivalent
to $T$ for all $\alpha$ automorphism of $\BBB^n$. Consequently,
Theorem~\ref{th:last}, (iii) says
that \emph{the standard multishift is homogeneous}. 


\bigskip
The authors thank J. Sarkar for useful discussions.


\begin{thebibliography}{xx}

\bibitem{Ar} W. Arveson: Subalgebras of $C\sp *$-algebras. III. Multivariable operator theory.  \emph{Acta Math.}  {\bf 181}  (1998), 159--228. 

\bibitem{Ar2} W. Arveson: The curvature of a Hilbert module over $\CCC[z_1,\cdots,z_d]$, \emph{ Proc. Natl. Acad. Sci. USA} {\bf 96}  (1999),  11096--11099.

\bibitem{BT2} Ch. Benhida, D. Timotin: Finite rank perturbations of contractions,  \emph{Integral Equations Operator Theory} {\bf 36}  (2000), 253--268. 

\bibitem{BES} T. Bhattacharyya, J. Eschmeier and J. Sarkar: Characteristic function
of a pure contractive tuple, \emph{Integral Equations  Operator Theory}, to appear.

\bibitem{ES} T. Bhattacharyya, J. Eschmeier and J. Sarkar: On completely non coisometric tuples and their
characteristic functions, preprint.


\bibitem{MS} G. Misra, N.S. Narsimha Sastry: Homogeneous tuples of operators and representations of some classical groups,\emph{ J. Operator Theory} {\bf 24}  (1990), 23--32.

\bibitem{Mu} V. M\"uller: \emph{Spectral Theory of Linear Operators}, Birkh\"auser Verlag,
Basel--Boston--Berlin, 2003

\bibitem{GP1} Gelu Popescu: Isometric dilations for infinite sequences of noncommuting operators,  \emph{Trans. Amer. Math. Soc.} {\bf 316 } (1989),  523--536. 

\bibitem{GP} Gelu Popescu: Characteristic functions for infinite sequences of noncommuting operators, \emph{ J. Operator Theory} {\bf 22}  (1989),  51--71. 

\bibitem{GP3} Gelu Popescu: Constrained multivariable operator theory, preprint, arXiv:math.OA/0507158.

\bibitem{RBD} B.V. Rajarama Bhat, T. Bhattacharyya, S. Dey:
Standard noncommuting and commuting dilations of commuting tuples,
\emph{Trans. Amer. Math. Soc.} {\bf 356} (2004), 1551--1568.

\bibitem{Ru} W. Rudin: \emph{Function Theory in the Unit Ball of $\CCC^n$}, Springer-Verlag, New York--Berlin, 1980.

\bibitem{SNF} B. Sz.-Nagy, C. Foias: \emph{Harmonic analysis of operators on Hilbert space}, North-Holland Publishing Co, 1970.

\end{thebibliography}
\end{document}